\documentclass[preprint,11pt]{elsarticle}

\usepackage{amssymb,amscd,latexsym}
\usepackage{amsthm}
\usepackage{amsmath}
\input xypic
\newcommand{\rar}{\rightarrow}
\newcommand{\lar}{\longrightarrow}
\newcommand{\llar}{-\kern-5pt-\kern-5pt\longrightarrow}
\newcommand{\surjects}{\twoheadrightarrow}
\newcommand{\injects}{\hookrightarrow}

\newtheorem{Theorem}{Theorem}[section]
\newtheorem{Lemma}[Theorem]{Lemma}
\newtheorem{Corollary}[Theorem]{Corollary}
\newtheorem{Proposition}[Theorem]{Proposition}
\newtheorem{Remark}[Theorem]{Remark}
\newtheorem{Example}[Theorem]{Example}

\newtheorem{Definition}[Theorem]{Definition}
\newtheorem{Question}[Theorem]{Question}
\def\sqr#1#2{{\vcenter{\hrule height.#2pt
        \hbox{\vrule width.#2pt height#1pt \kern#1pt
            \vrule width.#2pt}
        \hrule height.#2pt}}}
\def\phi{\varphi}
\def\demo{\noindent{\bf Proof. }}
\def\square{\mathchoice\sqr64\sqr64\sqr{4}3\sqr{3}3}
\def\qed{\hspace*{\fill} $\square$}

\def\xx{{\bf x}}

\def\TT{{\bf T}}

\def\XX{{\bf X}}

\def\ff{{\bf f}}

\def\ff{{\bf f}}

\def\hh{{\bf h}}

\def\hht{{\rm ht}\,}

\def\ker{{\rm ker}\,}

\newcommand{\ecod}{\mbox{\rm ecod }}

\def\rk{{\rm rank}\,}

\def\pp{{\mathbb P}}

\def\fd{{\cal D}}
\def\fp{\mathfrak{D}}
\def\fz{\mathfrak{Z}}

\journal{Journal of Algebra}

\begin{document}

\begin{frontmatter}



\title{Syzygies of differentials of forms}


\author{Isabel Bermejo\fnref{label1}}
\fntext[label1]{Partially supported by the Ministerio de Educaci\'on y
Ciencia - Espa\~na (MTM2010-20279-C02-02)}

\address{Facultad de Matem\'aticas, Universidad de La
Laguna, 38200 La Laguna, Tenerife, Canary Islands, Spain}


 \author{Philippe Gimenez\fnref{label2}}
\fntext[label2]{Partially supported by the Ministerio de Educaci\'on y
Ciencia - Espa\~na (MTM2010-20279-C02-02)}

 \address{Philippe Gimenez, Departamento de Algebra, Geometr\'{\i}a y
Topolog\'{\i}a, Facultad de Ciencias, Universidad de Valladolid,
47005 Valladolid, Spain}

\author{Aron Simis\fnref{label2}}
\fntext[label2]{Partially supported by the Ministerio de Educaci\'on y
Ciencia - Espa\~na (MTM2010-20279-C02-02); this author thanks the Universities of La Laguna and Valladolid for
hospitality and facilities provided during the preparation of this work}

\address{Departamento de Matem\'atica, Universidade
Federal de Pernambuco,\\ 50740-560 Recife, PE,
Brazil}



\begin{abstract}
Given a standard graded polynomial ring $R=k[x_1,\ldots,x_n]$ over a field $k$
of characteristic zero and a graded $k$-subalgebra $A=k[f_1,\ldots,f_m]\subset R$,
one relates the module $\Omega_{A/k}$ of K\"ahler $k$-differentials of $A$ to the transposed Jacobian module
$\mathcal{D}\subset \sum_{i=1}^n R dx_i$ of the forms $f_1,\ldots,f_m$ by means of a {\em Leibniz map}
$\Omega_{A/k}\rar \mathcal{D}$ whose kernel is the torsion of $\Omega_{A/k}$.
Letting $\fp$ denote the $R$-submodule generated by the (image of the) syzygy module of $\Omega_{A/k}$
and $\fz$ the syzygy module of $\mathcal{D}$, there is a natural inclusion $\fp\subset \fz$ coming
from the chain rule for composite derivatives.
The main goal is to give means to test when this inclusion is an equality -- in which case one says that
the forms $f_1,\ldots,f_m$ are {\em polarizable}. One surveys some classes of subalgebras  that are
generated by polarizable forms.
The problem has some curious connections with constructs of commutative algebra, such as the Jacobian ideal,
the conormal module and its torsion, homological dimension in $R$ and syzygies, complete intersections and
Koszul algebras.
Some of these connections trigger questions which have interest in their own.
\end{abstract}

\begin{keyword}
Syzygies\sep K\"ahler differentials\sep polarizability\sep Jacobian ideal\sep homological dimension\sep Polar map


\MSC 13C40\sep 13D02\sep 13B22\sep 13N05\sep 14M05

\end{keyword}

\end{frontmatter}





\section*{Introduction}

The subject envisaged  turns
out to be a special case of a general module
approximation situation, which we now state in a precise way.

We are given a field
$k$ of characteristic zero and an integral domain
$A$ of finite type over $k$ admitting an  embedding $A\simeq
k[\ff]=k[f_1,\ldots ,f_m]\subset R$, where $R=k[\xx]=k[x_1\ldots
,x_n]$ is a polynomial ring.
 Let $\fd$ be a finitely generated torsion free $R$-module
with an embedding $\fd\subset R^n\simeq A^n\otimes_AR$ that reflects the nature of the
above embedding $A\subset R$ - this vagueness is to be clarified in
each particular event. The question is whether one can
``approximate'' a free presentation of $\fd$ over $R$ by means of a
free presentation of a well-established $A$-module $\rm D$.

More precisely, let $0\rar \fz\rar R^m\rar \fd\rar 0$ be a
presentation of $\fd$ based on the generators of the embedding
$\fd\subset R^n$. One looks for an appropriate finitely generated
$A$-module $\rm D$ along with a presentation $0\rar {\rm Z}\rar A^m\rar {\rm D}\rar 0$
such that the $R$-submodule ${\mathfrak D}:={\rm Im}({\rm Z}\otimes_A R)\subset A^m\otimes_A R\simeq R^m$
generated by ${\rm Z}$ approximates $\fz\subset R^m$.

Besides looking for a sufficiently ubiquitous $\rm D$, one ought to
conceive a notion of ``approximation''. Natural features are:
\begin{enumerate}
\item[{\rm (i)}] (Inclusion) ${\mathfrak D}\subset \fz$;

\item[{\rm (ii)}] (Rank) $\rk_R({{\mathfrak D}})=\rk_R({\fz})$;

Further, since ${\fz}$ is a reflexive $R$-module, we may want that:

\item[{\rm (iii)}] (Depth) ${{\mathfrak D}}$ be a reflexive $R$-module;

Or else, one may wish that ${\fz}$ be closely approximated in low
codimension, such as:

\item[{\rm (iv)}] (Low codimension) ${{\mathfrak D}}={\fz}$ locally in codimension one
(i.e., $\hht ({\mathfrak D}:\fz)\geq 2$);

\item[{\rm (v)}] (Contraction) $A^m\cap \fz={\rm Z}$ in the natural inclusion of
$A$-modules $A^m\subset R^m$ induced by the ring extension $A\subset
R$.

\end{enumerate}

As is well-known, conditions (i) through (iii) imply that
${\mathfrak D}:_R\fz$ is height one unmixed or else ${\mathfrak D}:\fz=R$,
hence adding condition (iv)
yields the equality ${\mathfrak D}={\mathfrak{Z}}$. In particular, they
imply condition (v). On itself, (v) is a minimality condition on the
approximation.

In this paper we deal with the following situation:
one takes $\fd$ as the transposed Jacobian module $\fd(\ff)$ of
the forms $\ff$, while for $\rm D$ one takes the module $\Omega_{A/k}$ of K\"ahler
$k$-differentials of $A$ and then compare the respective modules of syzygies $\fz$ and $\fp$.
One should note that while $\Omega_{A/k}$  depends
only on the $k$-algebra $A$ and not on any particular presentation -- such as
$A\simeq k[\ff]\subset k[\xx]$ -- $\fd(\ff)$ depends on the choice
of the forms $\ff$.
What saves the  face of the transposed Jacobian module regarding
this instability is a {\em Leibniz map} $\Omega_{A/k} \rar \fd(\ff)$.
Although this map depends on $\ff$, its kernel is ultimately uniquely defined
and coincides with the $A$-torsion submodule of $\Omega_{A/k}$ (see  Theorem~\ref{Leibniz}).

We will say that $\ff$ is {\em polarizable} if ${\mathfrak D}=\fz$ in
the foregoing notation.

From a strict point of view, the problem considered in this paper has already been stated in \cite{jac, BGS}.
However, while in the latter one took a purely combinatorial approach, here one gets entangled in
several module theoretic queries.
As will be seen, this makes up for quite some difference between this and the previous references.

We focus  on a set of arbitrary forms of degree $2$.
Actually some of the preliminary results in this paper hold  for
forms $\ff$ of any fixed degree, but as we will see there is little hope to come around polarizability
in such generality.
From another end, polarizability is hardly a matter of
``general'' embeddings $A\subset R$. To illustrate this point,
consider sets of general
quadrics defining some of the well-known rational maps
$\pp^{n-1}\dasharrow \pp^{m-1}$: these  may fail to be polarizable even
when the image of such maps is a smooth variety (see
Example~\ref{smooth_quadrics}). Thus, there seems to be a strong
correlation to sparcity and non-genericity, perhaps nearly the same
way that birationality is thus related -- in fact, in \cite{BGS}
there are unexpected connections between these two notions and the
normality of $A$ in case $\ff$ are monomials of degree $2$. As we
guess, the lack of sparcity often imposes a high initial
degree on the defining equations of $A$, an early obstruction to
polarizability due to the relatively small initial degree of the
syzygy module $\fz$.

We next describe the contents of each section.

The first section contains the setup and definitions.
The main statement is a result bringing up
 the Jacobian ideal of $A$
as a tool to approximate the two modules.
We state a general result about modules of the same rank and
Fitting ideals (Lemma~\ref{conductors}).

In the second section we introduce the Leibniz map $\lambda:\Omega_{A/k} \rar \fd(\ff)$.
The main result (Theorem~\ref{Leibniz}) shows a tight relationship between the kernel
of $\lambda$, the torsion of $\Omega_{A/k}$ and the ``contraction'' of the differential
syzygies $\fz$ to $A$.
We then relate this result to polarizability, by showing that under a certain contractibility
hypothesis, polarizability implies the reflexiveness of the torsionfree $A$-module $P/P^{(2)}$,
where $A\simeq k[\TT]/P$.

The third section is dedicated to a discussion about when the algebra $A=k[\ff]$ is a complete intersection
on the presence of  polarizability.
As it turns out, this is the case nearly exactly when the transposed Jacobian module of $\ff$ has homological
dimension at most $1$.
There are some variations on this theme assuming conditions on $A$ of possibly unexpected nature.

In the next section we discuss the relationship between the ordinary syzygies of $\ff$ and their differential
syzygies, thereby introducing a framework where both fit into a basic exact sequence.
A consequence is the ability to give an alternative way of looking at the differential syzygies
as related to the ordinary syzygies.
As an application we show that by restricting the generating degrees of $P$, polarizability,
and hence, the complete intersection property, in fact follows from the assumption that the
homological dimension of $\fd(\ff)$ is at most $1$.

The last two sections are dedicated to a detailed analysis of selected examples with varied behavior,
and also examples of a more structured nature coming from constructs in algebraic geometry.
At the end we state some open questions whose answers look essential for further development
of the subject.

\section{Preliminaries}

We establish the basic setup, drawing upon the
notation stated in the introduction.
Thus, $A=k[\ff]\subset R$ denotes  a finitely generated $k$-subalgebra of the polynomial
ring $R$ and  $A\simeq
k[\TT]/P$  by mapping $T_j\mapsto f_j$.

Recall the well-known conormal sequence
\begin{equation*}
P/P^2\stackrel{\delta}{\lar}\sum_{j=1}^m A\, dT_j\lar
\Omega_{A/k}\rar 0,
\end{equation*}
where $\delta$ maps a generator of $P$ to its differential modulo $P$.
More exactly, upon choosing a generating set of $P$, the image of $\delta$ is the
$A$-submodule generated by the image of the transposed Jacobian
matrix over $k[\TT]$ of a generating set of $P$.
Since $P$ is a prime ideal and char$(k)=0$, we know the kernel of the leftmost map
above. Therefore, we will focus on the exact sequence
\begin{equation}\label{conormal_sequence}
0\rar P/P^{(2)}\stackrel{\delta}{\lar}\sum_{j=1}^m A\, dT_j\lar
\Omega_{A/k}\rar 0,
\end{equation}
where $P^{(2)}$ stands for the second symbolic power of $P$.

Throughout, we let $\fp\subset \sum_{j=1}^m R\,
dT_j$, denote the $R$-submodule generated by the image $\delta(P/P^{(2)})$
through the embedding $\sum_{j=1}^m A\, dT_j\subset
\sum_{i=1}^m R\, dT_j$ induced by the inclusion $A\subset R$.
Then $\fp$ is generated by the vectors $\sum_j \frac{\partial F}
{\partial T_j} \,(\ff)\,dT_j$, where $F$ runs through a set of
generators of $P$. On the other hand, by the chain rule of
composite derivatives, if $F\in P$ then $\sum_{j=1}^m \frac{\partial
F}{\partial T_j} (\ff)\, df_j=0$. This means that $\fp\subset \fz$,
where $\fz$ is the first syzygy module of the
differentials $d\ff$. In other words, the elements of $\fp$ are
relations of the transposed Jacobian matrix of $\ff$.
These ideas have also been discussed in a different context in \cite[Section 1.1]{Mich}
and even earlier in \cite[Main Lemma 2.3 (i)]{gauss}.
Borrowing from this line of thought, the present problem asks when a natural short complex
involving differentials is exact.

\begin{Definition}\rm  The elements of $\fp$ will be called
{\it polar syzygies\/} of $\ff$ (or of the embedding $A\subset R$ if
no confusion arises), while $\fp$ itself is named the {\it polar
syzygy module\/} of $\ff$. For the sake of comparison, we call $\fz$
the {\it differential syzygy module\/} of $\ff$ - though this is
actually the syzygy module of the differentials $d\ff$.
 Accordingly, we say that $\ff$ (or the embedding $A\subset R$) is
 {\it polarizable\/} if $\fp=\fz$.
\end{Definition}
The definition adopted here is the same as in \cite{BGS} but
slightly departs from the original definition given in \cite{jac}.
We will have a chance to deal
with the slight difference later in this work -- see also the
comment after Example~\ref{generators_drop}.
However,  in general it remains a bewildering piece.

\begin{Example}\label{failure_of_invariance}\rm
To see how the concept may actually depend on the
choice of the generators $\ff$, consider the following non-extravagant example:
$A=k[x^2, x^3]\subset k[x]$. Here, the differentials are $2xdx$ and
$3x^2dx$, respectively, hence the syzygy module of this set of
generators is the cyclic module generated by the single vector
$(3x,\, -2)^t$, while the polar module $\fp$ is generated by the
vector $(3x^4,\, -2x^3)^t$ that lies deep inside $\fz$. Thus, there
is an inclusion as it should be, but not an equality. On the other
hand, the transposed Jacobian module $\fd(\ff)$ is minimally
generated by the differential $2xdx$, hence, for this generator,
$\fz=\{0\}$. But, of course, $A$ itself cannot be generated over $k$
by less than $2$ elements, so we still have the same $\fp\neq \{0\}$
as before! Thus the theory depends on fixing the set $\ff$ of
generators of the $k$-algebra $A$, fixing $\fz$ as the kernel of the
presentation of $\fd(\ff)$ on the generators $d\ff$ and, likewise
fixing $\fp$ as the $R$-submodule generated by those polar syzygies
obtained by evaluating a set of generators of $P$ which are polynomial
relations of $d\ff$.
\end{Example}

Of course, some of the particular features of the above example
will not be present if one sticks to the truly {\em unirational} case,
i.e., when $A$ is generated by forms of a fixed degree.
But even then  polarizability may fail to be
an invariant property of the $k$-isomorphism class of $A$ (see
\cite[5.5]{BGS}).
Nonetheless a definite advantage of dealing with
homogeneous generators, not necessarily of the same degree,
is that, as a consequence of the Euler map, their $k$-linear independence is tantamount to that of the
corresponding differentials (see the setup in the proof of Proposition~\ref{EulerJacobiKoszul}).

\medskip

A trivial example of a polarizable embedding is the case when $A$
itself is a polynomial ring over $k$, i.e., when $\ff$ are algebraically
independent over $k$. This is scarcely of any interest for
the problem and we will rather look at genuine cases of the
notion.
For these, an important role will be played by the {\it Jacobian ideal\/} ${\mathfrak
J}={\mathfrak J}_A$ of $A$, which is read off the presentation
$A\simeq k[T_1,\ldots,T_m]/P$ by taking the ideal on $A$ generated
by the $g$-minors of the Jacobian matrix of a set of generators of
$P$, where $g=\hht P$. It is well-known that ${\mathfrak J}$
coincides with the $d$th Fitting ideal of the module of K\"ahler
$k$-differentials, where $d=\dim A$, hence is independent of the
presentation of $A$ and the choice of generators of $P$. We will
write ${\mathfrak J}R\subset R$ for the extended ideal in
$R=k[\xx]$.

\smallskip

In order to involve the Jacobian ideal we
will rely on a lemma of independent general interest, of
which we found no explicit mention, much less a proof, though it is conceivable that
it may be a piece of folklore.

To place it in a familiar framework, let $R$ be an integral domain and let
$E$ be finitely generated $R$-module with rank $r$; then some power of
the $r$th Fitting ideal $I\subset R$ of $E$ annihilates its torsion $\tau(E)$.
This assertion follows simply from the fact that $\tau(E)_P=\{0\}$ for
every prime ideal $P\not\in V(I)$.
In particular, if $E\surjects F$ is a surjective homomorphism of $R$-modules
of the same rank $r$  then
some power of the Fitting ideal of $E$ annihilates $\ker(E\surjects F)$.

The next result shows that if, moreover, $E$ and $F$ are generated by sets of the same cardinality then
$I$ itself already annihilates the kernel.
We state the result in the following format, which is more appropriate to the subsequent development.

\begin{Lemma}\label{conductors}
Let $R$ be an integral domain and let $M\subset N\subset R^m$ be
finitely generated submodules of a free module, having the same rank
$g$. Let $I\subset R$ denote the Fitting ideal
of order $m-g$ of the cokernel $R^m/M$. Then
$I\subset M:N$.
\end{Lemma}
\demo  The proof consists of a simple application of Laplace rule
for computing determinants. Namely, note that $I$ can be taken to be
the ideal generated by the $g\times g$ minors of the matrix whose
columns are the generators of $M$ expressed as linear combinations
of the canonical basis of $R^m$. Thus, let $\Delta\in I$ denote a
nonzero determinant thereof. We may assume for simplicity that it is
the determinant of the $g\times g$ submatrix on the upper left
corner. Given any $i=g+1,\ldots, m$, consider the following
$(g+1)\times g$ submatrix of the columns generating $M$:

$$\left(\begin{array}{cccc}
u_{1,1} & u_{1,2} & \cdots & u_{1,g}\\
u_{2,1} & u_{2,2} & \cdots & u_{2,g}\\
\vdots &\vdots  &  & \vdots \\
u_{g,1} & u_{g,2}  &\cdots & u_{g,g}\\
u_{i,1} & u_{i,2}  &\cdots & u_{i,g}
\end{array}
\right)
$$
Now, given any column generator of $N$, right border the above
matrix with the corresponding entries $v_1,\ldots,v_g,v_{i}$ of this
column to get a $(g+1)\times (g+1)$ matrix whose columns are
elements of $N$. These generate a submodule of $N$, hence has rank
at most $g$. Therefore the corresponding $(g+1)\times (g+1)$
determinant vanishes. Developing this determinant by Laplace along
the bottom row, one finds
\begin{equation}\label{expression}
\Delta v_{i}=\Delta_{\hat{1}2\ldots gv} u_{i,1}+
\Delta_{1\hat{2}3\ldots gv} u_{i,2}+\cdots
+\Delta_{12\ldots\hat{g}v} u_{i,g},
\end{equation}
 where
$\Delta_{12\ldots \hat{j}\ldots gv}$ denotes the $g$-minor obtained
by replacing the $j$th column with column $v$.

If, on the other hand, $i\in \{1,\ldots, g\}$ then we obtain again a
$(g+1)\times (g+1)$ matrix by first bordering the initial $g\times
g$ submatrix with the entries $v_1,\ldots,v_g$ and then repeating
the $i$th row of this matrix on the bottom. Clearly, this
determinant is zero; developing it as before along the repeated row,
we find a similar expression as (\ref{expression}), with the same
fixed $g$-minors as multipliers.

This shows that the entire column generator $v\in N$ is conducted by
$\Delta$ inside $M$.
\qed

\begin{Corollary}\label{jacobian_ideal_is_bona_fide}
Keeping the previous notation,  let
${\mathfrak J}\subset A$ stand for the Jacobian ideal of $A$. Then
${\mathfrak J}R\subset \fp:_R\fz$. Moreover, $\ff$ is polarizable if and only
if $\,\fp:_{R^m} {\mathfrak J}R=\fp$.
\end{Corollary}
\demo We have $\fp\subset \fz\subset \sum_{j=1}^m R\, dT_j$ and
$\rk_R(\fp)=\rk_A(P/P^{(2)})=\hht P =\rk_R(\fz)$ (see \cite[2.3]{BGS}.
By definition and the well-known fact that determinants commute with base change,
${\mathfrak J}R$ is the Fitting ideal of $R^m/\fp$ of order $m-\rk_R(\fp)$.
We apply Lemma~\ref{conductors} with $R^m=\sum_{j=1}^m R\, dT_j$, $M=\fp$, $N=\fz$,
where $I={\mathfrak J}R$. This proves the first statement.
The second assertion is a consequence of a general property of modules:
if $N$ is the kernel of a homomorphism $\eta:R^m\rar R^n$ of free
modules, then  $N:_{R^m} J=N$ for any
nonzero ideal $J\subset R$. Indeed, let $0\neq a\in J$ and $v\in R^m$ with $av\in N$.
Applying $\eta$ yields $a\eta(v)\in \eta(N)=\{0\}$, hence $\eta(v)=0$, i.e., $v\in N$,
as claimed.
Now, to get the second statement of the corollary, apply this fact with $N=\fz$.
 \qed

\bigskip

Of course, in the statement of the corollary one could replace
${\mathfrak J}R$ by any nonzero ideal ${\mathfrak a}\subset R$ such
that ${\mathfrak a}\subset \fp:_R\fz$. The point is that the Jacobian
ideal is a {\em bona fide} test ideal for polarizability that avoids
knowing $\fz$ a priori.  In this vein, it
would be interesting to know a priori when $\fz=(\fp)^{**}$
(bidual). Of course, it is often the case that $\fp$ is reflexive
but $\fp\neq\fz$ locally in codimension one.

For further insight into the role of the Jacobian ideal vis-\`a-vis polarizability
we refer to the  next section, where the role of the K\"ahler differential forms
is emphasized.

\section{The Leibniz map}

Recall the $R$-module $\mathfrak{Z}=\ker(\sum_{j=1}^m R\, dT_j\surjects \fd(\ff))$
and the conormal exact sequence of the module of K\"ahler $k$-differentials of $A$
as in (\ref{conormal_sequence}):
$$0  \rar  P/P^{(2)}  \stackrel{\delta}{\lar}  \sum_{j=1}^m A\,dT_j   \stackrel{\pi}{\lar}   \Omega_{A/k}
\rar  0.$$

\begin{Theorem}\label{Leibniz}
There is an $A$-module homomorphism $\lambda:\Omega_{A/k}\rar \fd(\ff)$
such that
$$\pi^{-1}(\ker(\lambda))=\pi^{-1}(\tau_A(\Omega_{A/k}))=(\mathfrak{Z})^c,$$
where $\tau_A(\Omega_{A/k})$ denotes $A$-torsion submodule of $\Omega_{A/k}$
and $(\mathfrak{Z})^c$ denotes the contraction of $\mathfrak{Z}$ via the natural inclusion
$\sum_{j=1}^m A\, dT_j\subset \sum_{j=1}^m R\, dT_j$ induced by $A\subset R$.
\end{Theorem}
\demo
We define the $A$-module homomorphism
$$\tilde{\lambda}:\sum_{j=1}^m A\, dT_j \surjects \sum_{j=1}^m A \,df_j
\subset\sum_{j=1}^m R \,df_j=\fd(\ff)\subset \sum_{i=1}^n R \,dx_i$$
by the association $\tilde{\lambda}(dT_j)=df_j$.
Then the chain rule for composite derivatives shows that $\tilde{\lambda}$ induces
an $A$-map $\lambda:\Omega_{A/k}\rar \fd(\ff)$.
Note that the image of $\lambda$ generates $\fd(\ff)$, although the map depends on $\ff$ -- i.e., on the
chosen embedding $A\subset R$, not just on $A$.
Thus, $\Omega_{A/k}$ is the only module in sight that depends only on $A$.

Now, $\lambda$ fits in a commutative diagram
of $A$-modules and $A$-homomorphisms
$$\begin{array}{ccccccccc}
0 & \rar & P/P^{(2)} & \stackrel{\delta}{\lar} & \sum_{j=1}^m A\,dT_j & \stackrel{\pi}{\lar} & \Omega_{A/k}
&\rar & 0\\[5pt]
&& \bigcap && \bigcap && \downarrow\lambda &&\\
 0 &\rar & \mathfrak{Z} & \lar & \sum_{j=1}^m R\, dT_j & \lar & \fd(\ff) &\rar
& 0
\end{array}.
$$
From this diagram one readily sees that $(\mathfrak{Z})^c\subset \pi^{-1}(\ker(\lambda))$.
On the other hand, one clearly has a surjection $\Omega_{A/k}/\tau_A(\Omega_{A/k})\surjects \Omega_{A/k}/\ker(\lambda)$
since $\fd(\ff)$ is $R$-torsionfree, hence also $A$-torsionfree.
Let $K$ (respectively, $L$) denote the fraction field of $A$ (respectively, of $R$).
Then
\begin{eqnarray*}
\rk_A(\Omega_{A/k}/\ker(\lambda))&=& \dim_K(\Omega_{A/k}/\ker(\lambda)\otimes_AK)=
\dim _L(R\,\Omega_{A/k}/\ker(\lambda)\otimes_RL)\\
&=&\dim _L ((\Omega_{A/k}/\ker(\lambda)\otimes_AK)\otimes_KL)\\
&=&\rk_R(\fd(\ff))=\dim A.
\end{eqnarray*}
It follows that the kernel of the surjection $\Omega_{A/k}/\tau_A(\Omega_{A/k})\surjects \Omega_{A/k}/\ker(\lambda)$
has rank $0$, hence must be the null module since $\Omega_{A/k}/\tau_A(\Omega_{A/k})$ is torsionfree.
Therefore, $\ker(\lambda)=\tau_A(\Omega_{A/k})$.

Now, again since $\Omega_{A/k}$ has rank $\dim A$, then $\rk_A(P/P^{(2)})=m-\dim A=\hht(P)$.
Therefore, the Fitting ideal of order $\dim A$ of $\Omega_{A/k}$ is the Jacobian ideal $\mathfrak{J}$ of $A$.
It follows that $\tau_A(\Omega_{A/k})=0:\mathfrak{J}^{\infty}$ (see \cite[Lemma 5.2]{elam99}).
Writing out this equality in terms of submodules of $\sum_{j=1}^m A\,dT_j$, yields the lifting
$\pi^{-1}(\tau_A(\Omega_{A/k}))=P/P^{(2)}:\mathfrak{J}^{\infty}$.
But the embedding $\sum_{j=1}^m A\,dT_j\subset \sum_{j=1}^m R\,dT_j$ induces an embedding
$P/P^{(2)}:\mathfrak{J}^{\infty}\subset \fp:(\mathfrak{J}R)^{\infty}$ which is preserved after contraction
back to $\sum_{j=1}^m A\,dT_j$.
On the other hand,
$$\fp:(\mathfrak{J}R)^{\infty}\subset \mathfrak{Z}:(\mathfrak{J}R)^{\infty}=\mathfrak{Z}, $$
since $\mathfrak{Z}$ is a second syzygy.
It follows that $\pi^{-1}(\tau_A(\Omega_{A/k}))\subset (\sum_{j=1}^m A\,dT_j)\cap \mathfrak{Z}=(\mathfrak{Z})^c.$
Collecting the pieces, we have
$$\pi^{-1}(\tau_A(\Omega_{A/k}))\subset (\mathfrak{Z})^c\subset \pi^{-1}(\ker(\lambda))= \pi^{-1}(\tau_A(\Omega_{A/k})),$$
thus proving the statement.
\qed
\begin{Definition}\rm
We call the map $\lambda$ the {\em Leibniz map}.
\end{Definition}

The following result is an immediate consequence of Theorem~\ref{Leibniz}.

\begin{Corollary}\label{easyfacts}
With the above notation, the following conditions are equivalent:
\begin{itemize}
\item[{\rm (i)}]
$P/P^{(2)}$ is contracted from $\fz$.
\item[{\rm (ii)}]
The Leibniz map $\lambda$ is injective.
\item[{\rm (iii)}]
$\Omega_{A/k}$ is torsion free.
\end{itemize}
Moreover, any of these conditions implies that $P/P^{(2)}$ is a reflexive $A$-module.
In particular, if $P/P^{(2)}$ is the contraction of its $R$-extension
$\fp$ and if $\ff$ is polarizable then $P/P^{(2)}$ is a reflexive
$A$-module.
\end{Corollary}

Additional applications will be given in the next section.

\begin{Remark}\rm
It is illuminating to compare the result in the above corollary to the one of
Corollary~\ref{jacobian_ideal_is_bona_fide}.
The reason why torsionfreeness of $\Omega_{A/k}$ does not imply polarizability
is that the former condition means that $P/P^{(2)}:_A \mathfrak{J}=P/P^{(2)}$,
and hence $(P/P^{(2)}:_A \mathfrak{J})\,R=\fp$,
while polarizability says that $\fp:_R \mathfrak{J}R=\fp$. Clearly, in general there is an
inclusion $(P/P^{(2)}:_A \mathfrak{J})\,R\subset (P/P^{(2)})\,R:_R \mathfrak{J}\,R=
\fp:_R \mathfrak{J}R$; thus, knowing that $\Omega_{A/k}$ is torsionfree  does not teach us
a lot more than the inclusion $\fp\subset \fp:_R \mathfrak{J}R$.
\end{Remark}
For explicit examples comparing the two properties, see (\ref{twisted_cubic}); and, in degree $2$, (\ref{old_jac}).

\section{Complete intersections}

In this part we focus on the case where $A$ is a complete intersection.
The first result is a criterion, in terms of polarizability,  for a normal almost complete intersection $A$ with ``small'' invariants
to be a complete intersection.

\begin{Proposition}\label{aci_perfect_cod2}
With the above notation, suppose that:
\begin{enumerate}
\item[{\rm (a)}] $A$ is a normal almost complete intersection{\rm ;}
\item[{\rm (b)}] $\dim A=\ecod A=2$, where $\ecod$ denotes embedding codimension.
\end{enumerate}
If the embedding $A\subset R$ is polarizable and $P/P^{(2)}$ is the contraction of its $R$-extension
then $A$ is a complete intersection in this embedding {\rm (}i.e., $P$ is generated by $2$
elements{\rm )}.
\end{Proposition}
\demo Let $A\simeq k[\TT]/P$ as before. Since $\dim A=2$ and $A$ is normal
(hence, satisfies $(S_2)$), then $A$
is Cohen--Macaulay. Since $A$ is an almost complete intersection
then $P^2=P^{(2)}$ (see \cite[(4.4)]{ram3}).
Then, by Corollary~\ref{easyfacts}, $P/P^2=P/P^{(2)}$ is reflexive, hence
Cohen-Macaulay because $\dim A=2$. By \cite[2.4]{He78}, $A$ is a
complete intersection. \qed

\begin{Question}\rm
In general the assumption that $P$ have deviation at most $1$
is essential for triggering the torsionfreeness of $P/P^2$
even if the conditions of (b) hold.
However, in the present case, $P$ is a particular prime ideal, so one asks
whether in the present context  the assumption of being almost complete intersection is
superfluous.
\end{Question}

\begin{Example}\label{twisted_cubic}\rm Let $A=k[x^3,x^2y,xy^2,y^3]\subset R=k[x,y]$ be the
parameters of the rational normal cubic. Then $A\simeq k[T_1,T_2,T_3,T_4]/P$, with $P=I_2(H)$,
where
$$
H=\left(\begin{array}{ccc}
T_1 & T_2 & T_3\\
T_2 & T_3 & T_4
\end{array}
\right).
$$
Therefore, $A$ is a normal Cohen-Macaulay almost complete intersection
of codimension $2$, hence $P/P^2$ is torsionfree, but not reflexive.
In fact, since $A$ is non-obstructed, $P/P^2$ is a proper reduction of
the double dual $(P/P^2)^{\ast\ast}$ (\cite[Corollary 4.5]{ram3}).
Also, $\Omega_{A/k}$ has nontrivial torsion submodule.
Note that $P/P^2$ is contracted from its $R$-extension.
Therefore, $A\subset R$ is not polarizable.
\end{Example}

\smallskip

Of a different flavor is the next result, naturally extending
\cite[Theorem 2.5, (a) $\Leftrightarrow$ (c)]{jac}.

\begin{Proposition}\label{homol_dim_one}
With the above notation, suppose that $\ff$ are homogeneous of the
same degree and that the embedding $A=k[\ff]\subset R$ is polarizable.
Writing $A\simeq k[\TT]/P$, assume that $\mu(P)=\mu({\fp})$.
The following conditions are equivalent:
\begin{enumerate}
\item[{\rm (a)}] The homological dimension of $\fd(\ff)$ is at
most one$\,${\rm ;}
\item[{\rm (b)}] $P$ is generated by a regular sequence {\rm (}i.e., $A$ is a
complete intersection{\rm )}.
\end{enumerate}
\end{Proposition}
\demo (a) $\Rightarrow$ (b) By a well-known stability preliminary,
the homological assumption means that $\fz$ is a projective module,
hence free. Since $\ff$ is assumed to be polarizable,  $\fp$ is free
as well (of rank $\hht P$ \cite[2.3]{BGS}). Now, by assumption
$\mu(P)=\mu({\fp})$. On the other hand, $\mu_A(P/P^2)=\mu(P)$
because $P$ is homogeneous. Then one gets $\mu_A(P/P^2)=\hht P=\rk _A(P/P^2)$.
It follows that $P/P^2$ is a free $A$-module and hence, $P$ is generated
by a regular sequence of forms by Ferrand--Vasconcelos theorem.

(b) $\Rightarrow$ (a) Since $P$ is generated by a regular sequence then
$P/P^2$ is a free $A$-module. The assumption $\mu(P)=\mu({\fp})$ then yields $\mu(\fp)=\hht P=
\rk _R(\fp)$. By polarizability, $\mu(\fz)=\rk_R(\fz)$, hence $\fz$ is free.
This shows that the homological dimension of $\fd(\ff)$ is at
most one.
 \qed

 \medskip

 The previous result stands on the assumption that $\mu(P)=\mu({\fp})$.
 Here is one situation where this equality holds.

\begin{Proposition}\label{no_generators_does_not_drop}
Suppose that $\ff$ are $k$-linearly independent forms of the same degree and, moreover,
that $P$ is generated in the same degree.
Then $\mu(P)=\mu({\fp})$.
\end{Proposition}
\demo
Morally, this comprises two steps. First,
$\mu(P)=\mu(P/P^{(2)})$, i.e., $P$ does not loose minimal generators
in the passage to the symbolic conormal module. Namely, we claim
that the kernel of the following map of $k$-vector spaces
$$ \frac{P}{(\TT)P}\surjects \frac{\kern3ptP/P^{(2)}}{(\TT)P/P^{(2)}}$$
vanishes. Indeed,  by the Zariski--Nagata differential criterion for
symbolic powers in characteristic zero, the elements of $P^{(2)}$
satisfy the property that all its partial derivatives belong to $P$.
Also, since $P^{(2)}$ is the $P$-primary component of the
homogeneous ideal $P^2$, it is homogeneous. Therefore,  the Euler
formula as applied to individual elements of $P^{(2)}$ yields
$P^{(2)}\subset (\TT)P$, as required.

Second step: $\mu(P)=\mu(\fp)$, i.e., $P$ -- or $P/P^{(2)}$, which
now amounts to the same -- does not loose minimal generators while
extending to an $R$-module. For this to hold we need the second assumption.
Thus, let $P$ be minimally generated by forms
$\{F_1,\ldots, F_r\}$ of the same degree. In any case $\fp$ is
generated by the evaluated differentials $\{dF_1(\ff),\ldots,
dF_r(\ff)\}$ as an $R$-module. But by the assumption on $\ff$ these
generators have the same standard $R$-degree, so the only way they
can fail to be minimal is that they be linearly dependent over the
base field $k$. However, such a dependence immediately implies that there are
$\lambda_1,\ldots,\lambda_r\in k$, not all zero, such that
$$\sum_{1\leq l\leq r}\,\lambda_l\, \frac{\partial F_l}{\partial
T_j}\in P, \; 1\leq j\leq m.$$ Since $\deg (\partial F_l/\partial
T_j)<\deg(F_l)$ for $\partial F_l/\partial T_j\neq 0$, we arrive at
a contradiction.
\qed

\medskip

The hypothesis on $P$ seems to be needed in general, as the following example
indicates.

\begin{Example}\label{generators_drop}\rm
Let $R=k[x_1,\ldots, x_6]$ and let
\begin{equation}\label{hex_triangle}
\ff=\{x_1x_2,\, x_2x_3,\, x_3x_4,\, x_4x_5,\, x_5x_6,\,
x_6x_1,\, x_1x_3,\, x_3x_5,\, x_5x_1\}.
\end{equation}
Then the presentation ideal of $k[\ff]$ is
a codimension $3$ almost complete intersection in which the quadrics
generate a maximal regular sequence, while the fourth generator
lives in degree $3$. It can be shown that $\ff$ is polarizable and
that $\fp$ is a free module generated by the differentials of the
three quadrics (see \cite[Example 5.21]{BGS}).
\end {Example}

It would be interesting to know, under the assumption that
$\ff$ are homogeneous of the same degree and perhaps also under
the assumption  of polarizability, when the equality $\mu(P)=\mu({\fp})$
holds. By a quirk this is the case if $\ff$ happens to be the set of
degree $2$ monomials corresponding to a connected bipartite graph -- see
\cite[Theorem 2.3]{jac}, where this falls in an indirect way from the main
result.

\begin{Remark}\rm To see how subtle the problem is, one can take the following bipartite
graph, whose edge-ideal is ideal theoretically entirely analogous to the above:
\begin{equation}\label{hex_squares}
\ff=\{x_1x_2,\, x_2x_3,\, x_3x_4,\, x_4x_5,\, x_5x_6,\,
x_6x_1,\, x_1x_7,\, x_3x_7,\, x_5x_7\}.
\end{equation}
Both this and (\ref{hex_triangle}) are Cohen--Macaulay ideals of
codimension $3$, with same graded Betti numbers.
However, here $\mu(P)=\mu(\fp)$ (by {\em loc. cit.} or
by direct computation).
\end{Remark}

We observe that in both (\ref{hex_triangle}) and (\ref{hex_squares}) the ideal $P$ is the
homogeneous defining ideal of an arithmetically normal projective variety.
Thus, in both cases one can deduce polarizability  from \cite[Theorem 5.10]{BGS}.

\medskip

We have the following immediate consequence.

\begin{Corollary}\label{homol_dim_one_and_same_degree}
Suppose that $\ff$ are homogeneous of the same degree and that the
embedding $k[\ff]\subset R$ is polarizable. If $P$ is generated in
the same degree and the homological dimension of $\fd(\ff)$ is at
most one then $A$ is a complete intersection.
\end{Corollary}
\demo One applies Proposition~\ref{homol_dim_one} along with
Proposition~\ref{no_generators_does_not_drop}.
\qed

\bigskip

We emphasize the role of polarizability in the last corollary. Thus,
e.g., if $\ff$ are the parameters defining the twisted
cubic in $\pp^3$, the hypotheses in the above statement are satisfied,
nevertheless $A$ is not a complete intersection, and indeed it is
not polarizable (see Example~\ref{twisted_cubic}).

At the other end of the spectrum, even  under the hypothesis of
polarizability the result is false in higher homological dimension.

\begin{Example}\rm
Let $R=k[x_1,x_2, x_3]$ and let $A=R^{(2)}\subset R$ be the
$2$-Veronese of $R$. Then this embedding is polarizable by
\cite[5.1]{BGS} and $A$ is defined by quadrics, while the projective
dimension of $\fd(\ff)$ is $2$.
\end {Example}

The question naturally arises as to when the two conditions in
Corollary~\ref{homol_dim_one_and_same_degree} actually imply that
$k[\ff]\subset R$ is polarizable (hence a complete intersection). In
the next section we
give an answer to this question under an additional condition (cf.
Proposition~\ref{homol_dim_one_and_quadrics}).

\section{The Euler--Jacobi--Koszul exact sequence}

In this section we still assume that $\ff$ are forms, but not necessarily of
the same degree.

We wish to relate more closely the polar and differential syzygies
of $\ff$ to its ordinary syzygies. Henceforth, for any set $\hh\subset R$, $Z(\hh)$ denotes the
first syzygy module of $\hh$.
Set $d_j=\deg(f_j), \,1\leq j\leq m$ and  $\tilde{\ff}=\{d_1f_1,\ldots, d_mf_m\}$.
Since char$(k)=0$, clearly $(\ff)R=(\tilde{\ff})R$.

\begin{Proposition}\label{EulerJacobiKoszul}
There is an exact sequence of $R$-modules
\begin{equation}\label{EJK}
0\rar {\mathfrak{Z}}\lar Z(\tilde{\ff})\stackrel{\mathfrak j}{\lar} \fd(\ff)\cap
Z({\xx}) \rar 0,
\end{equation}
where ${\mathfrak j}$ sends a syzygy $(g_1,\ldots, g_m)$ of $\tilde{\ff}$ to the
element $\sum_{j=1}^m g_j df_j$.
In particular, if $\ff$ is polarizable then there is an exact
sequence of $R$-modules
\begin{equation}\label{polarsequence}
0\rar \fp\lar Z(\tilde{\ff})\lar \fd(\ff)\cap Z(\xx) \rar 0.
\end{equation}
\end{Proposition}
\demo  The restriction of the
Euler map $\epsilon: \sum_{i=1}^n R\,dx_i\rar R$, where
$\epsilon(dx_i)=x_i$, induces an exact sequence of $R$-modules
$$0\rar \fd(\ff)\cap Z(\xx)\rar \fd(\ff)=\sum_{j=1}^m R\,df_j \rar (\tilde{\ff})R=(\ff)R\rar 0.$$
This in turn fits in a snake diagram:

{\small
$$\begin{array}{ccccccccc}
&&  &&  && 0 &&\\[5pt]
&&  &&  && \downarrow &&\\[5pt]
&& 0 && 0 && \fd(\ff)\cap Z({\xx}) &&\\[5pt]
&& \downarrow && \downarrow && \downarrow&&\\[5pt]
0 & \rar & {\mathfrak{Z}} & {\lar} & \sum_{i=1}^m R\, dT_j &
\stackrel{dT_j\mapsto df_j}{\lar} & \fd(\ff)
&\rar & 0\\[5pt]
&& \downarrow && \Vert && \downarrow\epsilon &&\\[5pt]
0 & \rar & Z(\tilde{\ff}) & {\lar} & \sum_{i=1}^m R\, dT_j & \stackrel{dT_j\mapsto d_jf_j}{\lar}
& (\ff)R &\rar & 0\\[5pt]
&& \downarrow && \downarrow && \downarrow&&\\[5pt]
&& Z(\tilde{\ff})/{\mathfrak{Z}} && 0 && 0 &&\\[5pt]
&& \downarrow &&  && &&\\[5pt]
&& 0 &&  && &&
\end{array},
$$
}
Then ${\mathfrak j}$ is the lifting to $Z(\tilde{\ff})$ of the inverse of
the connecting isomorphism in the kernel-cokernel sequence by the snake lemma.
To make this map explicit, note that an element $\sum_{j=1}^m g_j df_j\in  \fd(\ff)\cap
Z({\xx})$ is characterized by the equation $\mathbf{g}\cdot\Theta(\ff)\cdot \xx^t=0$
or by its transpose, where $\Theta$ denotes the Jacobian matrix of $\ff$.
This readily gives the way ${\mathfrak j}$ acts.
\qed

\begin{Definition}\rm
The exact sequence (\ref{EJK}) could be called the {\it
Euler--Jacobi--Koszul syzygy sequence\/}, while $\fd(\ff)\cap Z({\xx})$ could accordingly
be dubbed the {\it module of Euler--Jacobian syzygies\/} of $\ff$
and ${\mathfrak j}$ the {\it Jacobi
map\/} of $\ff$.
\end{Definition}

Both exact sequences are ways of showing, in characteristic zero,
that the syzygies of a polarizable set of forms of the same degree has
a fixed structure in differential terms.

\bigskip

Here is a computational view of the above syzygy sequence in terms
of the involved matrices, emphasizing the Jacobi map. For the sake of simplicity we assume that $\ff$ are
forms of the same degree and for the sake of lighter reading we
write $R^n=\sum_{i=1}^m R\, dx_i$, $R^m=\sum_{j=1}^m R\, dT_j$.

Denoting by $\Theta^t$ the transposed Jacobian matrix of $\ff$, one has maps

$$\begin{array}{ccccccccc}
\kern10pt R^q \kern-6pt&&&&&&&&\\
 & \kern-10pt\psi{\searrow} &&&&&&\\
&&  \kern-15pt R^p & \kern-35pt\stackrel{\phi}{\lar} & \kern-15pt
R^m \lar & (\ff)
& \rar & 0\\[10pt]
&& \Theta^t\cdot\phi{\searrow} &\kern5pt\swarrow\kern5pt\Theta^t&&&&\\
&&& \kern-40pt R^n &&&&&
\end{array}
$$
where $\phi$ denotes the syzygy matrix of $\ff$ and $\psi$ denotes the syzygy matrix
of $\Theta^t\cdot\phi$.
Then
$${\rm im}(\Theta^t\cdot\phi)\subset {\rm im}(\Theta^t)\cap
{\rm im}(\bigwedge^2R^n\stackrel{\kappa}{\lar} R^n),$$ with $\kappa$
denoting the first map of the Koszul complex on $\xx$. Thus, the columns of the
matrix $\Theta^t\cdot\phi$ generate the image of the map ${\mathfrak
j}$. Viewed this way,  $\fz$ is also the submodule of
$R^m$ generated by the columns of the product matrix
$\phi\cdot\psi$, though not minimally -- in fact, quite often some columns
may be null.

Certainly, this is not how one would like to compute generators of $\fz$, as it depends on
computing syzygies of an even more involved module.
Its purpose is mostly theoretical, indicating the intermediation of the syzygies of $\ff$, thus
leading us to make assumptions on those syzygies.

\medskip

Let indeg$\,_T(E)$ the initial degree of a graded module over a graded ring $T$.
One has the easy

\begin{Lemma}\label{obviety}
Let $\ff\subset R$ be $k$-linearly independent forms of the same degree $\geq 2$.
Then {\em indeg}$\,_R(\fp)\geq 2$ and {\em indeg}$\,_R(\fz)\geq 2$.
Moreover, suppose that $\fp$ is generated in fixed degree $d$, that {\em indeg}$\,_R(\fz)\geq d$ and that
$\mu(\fp)\geq \mu(\fz)$.
Then $\ff$ is polarizable.
\end{Lemma}
\demo
Since $\ff$ are $k$-linearly independent forms of the same degree, one has an isomorphism of graded $k$-algebras
 $k[\TT]/P\simeq A=k[\ff]\subset R$, with
$P$ is a homogeneous ideal generated in
degree $\geq 2$. Clearly then the image of $P/P^2$ in $\sum_j A\,dT_j$ is a graded $A$-module
generated in degree $\geq 1$. Recall that $\fp$ is generated by these generators further evaluated on $\ff\subset R$.
Since $\deg(\ff)\geq 2$, we clearly have indeg$\,_R(\fp)\geq 2$.

As for $\fz$, as seen above the entire syzygy module
$\fz$ is generated by the columns of the product matrix
$\phi\cdot\psi$, where $\phi$ is the syzyzgy matrix of $\ff$, hence its entries have $R$-degree at least $1$.
Taking a minimal set of generators of Im$(\Theta^t\cdot\phi)$,
the entries of the syzygy matrix $\psi$ of
$\Theta^t\cdot\phi$ will be forms of positive degree. It follows that
the columns of $\phi\cdot\psi$ have degree at least $2$, hence $\fz$
is generated in degree at least $2$.

The last assertion is clear since then all the minimal generators of
$\fp$ must be minimal generators of $\fz$ and the latter can thereby have no
other minimal generators.
\qed

\smallskip

An example of application of this sort of ideas is as follows.

\begin{Proposition}\label{homol_dim_one_and_quadrics}
With the above notation, suppose that $\ff$ are $k$-linearly
independent forms of degree $2$ and set $k[\ff]\simeq k[\TT]/P$ as
before. Assume that:
\begin{enumerate}
\item[{\rm (a)}] $P$ is minimally generated by forms of degree $2\,${\rm ;}
\item[{\rm (b)}] The homological
dimension of $\fd(\ff)$ is at most one.
\end{enumerate}
Then $A$ is a polarizable complete intersection.
\end{Proposition}
\demo By (b)  $\fz$ is a free module, i.e.,
$\mu(\fz)=\rk \,\fz$. Since $\rk \,\fp=\rk \,\fz$, we have
$\mu(\fp)\geq \mu(\fz)$. By Lemma~\ref{obviety}, $\fp=\fz$ must be
the case. That $A$ is besides a complete intersection follows from
Corollary~\ref{homol_dim_one_and_same_degree}. \qed

\begin{Remark}\rm Assumption (a) in the previous proposition is essential -- see
Example~\ref{nonpolarizable_hypersurface} below.
This shows that, in general, normal complete intersections are not polarizable.
\end{Remark}

\section{Illustrative examples}

The following is a selection of  examples to help visualize the
theory.

\begin{Example}\label{degree_higher_than_two}\rm If the degrees of the
generators $\ff$ of $A$ are higher than $2$, polarizability is a
rare phenomenon even if $A$ is a homogeneous isolated singularity.
To see this, consider a homogeneous version of
Example~\ref{failure_of_invariance}: $A=k[x_1^3,
x_1^2x_2,x_2^3]\subset k[x_1,x_2]$. A simple calculation shows that $\fz$ is (cyclic)
generated in $R$-degree $3$, while $\fp$ is (cyclic) generated in
$R$-degree $6$ because the presentation ideal of $A$ over $k$ is
generated by a form of degree $3$ in the presentation variables.
Perhaps more conceptual is the case of the parameters of the twisted cubic as shown in
Example~\ref{twisted_cubic}.
\end{Example}

\begin{Example}\label{smooth_quadrics}\rm If $\ff$ are general forms of degree
$2$ then  polarizability may fail. Indeed, if $\ff$ are $5$ general quadrics in
$k[x_0,x_1,x_2]$ then $k[\ff]$ is (up to degree renormalization) the
homogeneous coordinate ring of a general projection $V$ to $\pp^4$
of the $2$-Veronese embedding of $\pp^2$ in $\pp^5$. It is known and
classical that $V$ is a smooth variety cut out by cubics, hence the polar
syzygy  module $\fp$ is minimally generated in degree $4$. But the
differential syzygy module $\fz$ is generated in degree $3$.
\end{Example}
This example teaches us quite a bit:

\smallskip

(1) $A$ is an isolated singularity (i.e., Proj$(A)$ is smooth), but it is not normal
nor Cohen--Macaulay as ${\rm depth} (k[\ff])=1$;

(2) $\ff$ are the Pfaffians of a skew-symmetric $5\times 5$ matrix, whose entries are then
necessarily linear forms. Therefore, the ideal $(\ff)$ is linearly presented.
This is in direct contrast with the situation where $\ff$ are monomials, in which case $\ff$ is always
polarizable (see \cite[Proposition 5.18]{BGS}).

In any case, one may observe that through the Jacobi map of Proposition~\ref{EulerJacobiKoszul} the linear syzygies of
$\ff$ correspond to the degree $2$ component of the kernel of the Euler map $\mathcal{D}(\ff)\rar (\ff)$.
Thus, if $\ff$ linearly presented then this kernel is generated in degree $2$,  a condition
which would be interesting to understand.
Unfortunately, this condition does not seem to have immediate impact on the degrees of $\fz$
as the generators of the latter will be non-minimal syzygies of $\ff$.

\begin{Example}\label{nonpolarizable_hypersurface}\rm (See \cite[Theorem 5.10]{BGS})
A normal hypersurface is not polarizable in general.
 Let $$A=k[x_1^2,x_1x_2,x_2x_3,x_3x_4,x_2x_4]\subset
k[x_1,x_2,x_3,x_4].$$ Then $A$ is a normal hypersurface of degree
$3$ with defining equation $T_2^2T_4-T_1T_3T_5$. Obviously, $\fp$ is
free (of rank one). $A$ is not polarizable since the transposed
Jacobian matrix of its generators has a nonzero relation in degree
$3$ while $\fp$ is generated in degree $4$. Actually, $\fz$ is
cyclic, hence free. Thus, though $\fp$ is isomorphic to $\fz$ as
abstract $R$-modules, it is not a second syzygy in its natural
embedding (the inclusion $\fp\subset \fz$ is not an equality in
codimension one exactly along the prime $(x_2)$).

A computation  shows that $P/P^2$ is the
contraction of $\fz$, hence is contracted from its $R$-extension as
well. Thus, contractibility and reflexivity together do not imply
polarizability, hence the last assertion in
Corollary~\ref{easyfacts} admits no weak converse. One may wonder
whether contractibility is always the case for a normal complete
intersection parameterized by $2$-forms. Note, however, that there may exist ideals in $A$
which are not contracted from their extensions in $R$, as is here
the case -- e.g., $I=(x_3x_4)A$ is not contracted from its extension
since $x_2^2\cdot x_3x_4\in IR\cap A\setminus A$.
\end{Example}

\begin{Example}\label{old_jac}\rm (\cite[p. 992]{jac})  $A$ is now generated by square-free
monomials and $\fp$ is not a reflexive $R$-module. Namely, one
takes
$$A=k[x_1x_2,x_2x_3,x_3x_4,x_4x_5,x_5x_6,x_6x_1,x_2x_4,x_2x_6]\subset
k[x_1,x_2,x_3,x_4,x_5,x_6].$$ Here one can show that $A$ is normal,
hence Cohen--Macaulay, i.e., $P$ is a codimension two perfect ideal.
Moreover, it can easily be seen that $P$ is an almost complete intersection.
It follows that $P^2=P^{(2)}$ (see \cite[(4.4)]{ram3}).
On the other hand, $A$ is locally regular in codimension $2$, hence
$\Omega_{A/k}$ is torsionfree (see \cite[4.6]{ram3}).
It follows from Corollary~\ref{easyfacts} that $P/P^2$ is both
reflexive and contracted.
 It was pointed out in \cite{jac} that $A$ is not polarizable.
An additional calculation gives
$\hht(\fp:\fz)\geq 2$, hence the two modules coincide in codimension one.
Of course $\fp$ fails to be reflexive since $A$ is not polarizable.
\end{Example}

\begin{Example}\rm Next is an example where nearly everything goes wrong although the
$k$-generators $\ff$ of $A$ in the embedding $A\subset R$ have some
inner symmetry:
$$A=k[x_1^2,x_3^2,x_1x_2,x_2x_3,x_3x_4,x_1x_4]\subset k[x_1,x_2,x_3,x_4].$$
The presentation ideal $P$ of $A$ is a codimension two ideal
generated by one quadric and three cubics. Clearly, the quadric
forbids the generators to be minors.
Therefore, $A$ is not Cohen--Macaulay, hence not normal either
(although it it is locally regular in codimension one). Here,
$P^2=P^{(2)}$ and $P/P^2$ is not a reflexive $A$-module. $A$ is not polarizable,
hence $\fp$ is not a reflexive $R$-module either -- as a direct computational check, e.g.,
one can see that $\hht(\fp:\fz)\geq 2$.
\end{Example}

\section{Structured ``geometric'' classes}

In the previous section most examples had no structured nature.
In this part we collect a few  examples of a more structured geometric nature.

Previously known results are as follows:

\medskip

\noindent{\sc Veronese embeddings of order $2$}

It has been proved in \cite[Corollary 5.1]{BGS} that the edge algebra of a
complete graph with a loop at every vertex is polarizable.
This means that the $2$ Veronese embedding of $\pp^n$ is given by a polarizable
parametrization.

\medskip

\noindent{\sc Segre embedding and its coordinate projections}

It was proved in \cite[Corollary 5.2]{BGS} that the edge algebra of a connected bipartite
graph is polarizable. This means that the homogeneous coordinate rings of the Segre
embedding $\pp^n\times \pp^m\injects \pp^N$ and of all its ``connected''
coordinate projections  are polarizable.

\subsection{Grassmannian of lines and scrollar parameterizations}\label{Grassmannian}

In this part we prove two new results.

\begin{Proposition}\label{grassmannian} Let $\XX:=(x_{kl})\, (1\leq k\leq 2,\, 1\leq l\leq m\geq 3)$ be a $2\times m$
generic matrix over $k$. Let $R=k[x_{kl}]$ and
$A=k[\ff]\subset R$ where $\ff$ is the set of the $2\times 2$ minors
of $\XX$.
Then $\ff$ is polarizable.
\end{Proposition}
\demo It is classical and well-known that a defining ideal
$P\subset k[\TT]$ of $A$ over $k$ is generated by the Grassmann--Pl\"ucker
(quadratic) relations of these minors and such $P$ is minimally
generated by ${{m}\choose{4}}$ such relations quadratic $\TT$-forms, one for each $2\times
4$ submatrix of $\XX$.
By Proposition~\ref{no_generators_does_not_drop}, $\fp$ is minimally generated by ${{m}\choose{4}}$ elements. We will show that any element
in $\fz$ can be expressed as an $R$-combination of these ${{m}\choose{4}}$ polar syzygies, i.e., $\ff$ is polarizable.

For this we will proceed by induction on the number $m$ of columns (not on the number $2m$  of variables!).
The assertion is (vacuously) verified for $m=3$ since the minors are algebraically independent in this case and $\fp=\fz=0$.

Thus, assume that $m\geq 4$. To stay clear, we write $\XX(m)$ for $\XX$, $\ff(m)$ for $\ff$ and, accordingly, write $P(m)$
for $P$ and $\fz(m)$ for $\fz$.
Note that the equality ${{m}\choose{4}}={{m-1}\choose{4}}+{{m-1}\choose{3}}$
is the numerical counterpart of saying that $P(m)=(P(m-1), Q(m-1))$,
where $Q(m-1)$ denotes the set of Grassmann--Pl\"ucker relations obtained from all $2\times 4$ submatrices of
$\XX(m)$ involving the last column.

On the other hand, the transposed Jacobian matrix of $\ff(m)$ has the following shape
{\small
$$^t \Theta(\ff(m))=
\left(
\begin{array}{c@{\quad\vrule\quad}ccccc}
 & \hbox{$x_{2m}$} &0 &0 &\cdots &0\\
  & 0 & \hbox{$x_{2m}$} &0 &\cdots &0\\
\hbox{$^t \Theta(\ff(m-1))_1^{m-1}$} &0 &0 & \hbox{$x_{2m}$} &\cdots &0\\
 & \vdots &\vdots &\vdots &\cdots &\vdots\\
 & 0 & 0 &0 &\cdots & \hbox{$x_{2m}$}\\[5pt]
\multispan6\hrulefill\\[2pt]
\boldsymbol 0 & -x_{21} & -x_{22} & -x_{23} &\cdots & -x_{2m-1}\\[5pt]
\multispan6\hrulefill\\[2pt]
 & \hbox{$-x_{1m}$} &0 &0 &\cdots &0\\
 & 0 & \hbox{$-x_{1m}$} &0 &\cdots &0\\
\hbox{$^t \Theta(\ff(m-1))_{m+1}^{2m-1}$} &0 &0 & \hbox{$-x_{1m}$} &\cdots &0\\
 & \vdots &\vdots &\vdots &\cdots &\vdots\\
 & 0 & 0 &0 &\cdots & \hbox{$-x_{1m}$}\\[5pt]
\multispan6\hrulefill\\[2pt]
\boldsymbol 0 & x_{11} & x_{12} & x_{13} &\cdots & x_{1m-1}
\end{array}
\right),
$$
}
where the left vertical block, interspersed with two rows of zeros, is the transposed Jacobian matrix of $\ff(m-1)$.
By the inductive hypothesis, $\fz(m-1)$ is generated by ${{m-1}\choose{4}}$ polar syzygies. Let's denote by
$\widetilde{\fz(m-1)}$ the submodule of $\fz(m)$ minimally generated by these ${{m-1}\choose{4}}$ polar syzygies
each stacked over a column of
$m-1$ zeros. Observe that $\widetilde{\fz(m-1)}$ is the submodule of $\fz(m)$ generated by the polar syzygies
arising from the minimal generators of $P(m-1)$, and that any element in $\fz(m)$ not involving the
right vertical block above (i.e., whose last $m-1$ entries are zero) belongs to $\widetilde{\fz(m-1)}$.

Let $\mathbf{h}\in R^{m-1}$ denote the vector whose coordinates the last $m-1$ coordinates of an arbitrarily given
vector $\mathbf{s}\in \fz(m)$.
By the shape of the matrix, $\mathbf{h}$ is a syzygy of the  matrix formed with  the $m$-th and last rows of the matrix
right vertical block. Up to a permutation of these rows and a sign, this matrix is just $\XX(m-1)$.
We know the minimally generating syzygies of this generic matrix from the Buchsbaum--Rim complex.
To wit, for every $2\times 3$ submatrix of $\XX(m-1)$ with columns $1\leq i<j<k\leq m-1$, take its $2\times 2$ signed minors.
Then let ${\bf m}_{ijk}\in R^{m-1}$ denote the vector whose coordinates are these minors in places $i$, $j$ and $k$, and zero elsewhere.
The set of these vectors minimally generate the module of syzygies of $\XX(m-1)$.

Accordingly, write $\displaystyle{{\bf h}=\sum_{1\leq i<j<k\leq m-1}q_{ijk}{\bf m}_{ijk}}$ for suitable $q_{ijk}\in R$.

\smallskip

From the other end, choosing a $2\times 3$ submatrix of $\XX(m-1)$ with columns $i,j,k$ amounts to picking a
$2\times 4$ submatrix of $\XX(m)$ involving the last column. The latter gives rise to a Grassmann--Pl\"ucker
equation $\mathbf{p}_{ijk}$ inducing a polar syzygy whose nonzero coordinates in the last $m-1$ slots are precisely
the (signed) minors of the $2\times 3$ submatrix in places $i$, $j$ and $k$ as above.

It follows that the last $m-1$ coordinates of
$\displaystyle{{\bf s}\,-\sum_{1\leq i<j<k\leq m-1}q_{ijk}{\bf p}_{ijk}}$ are zero, hence this vector belongs to $\widetilde{\fz(m-1)}$.
Therefore, ${\bf s}$ is an $R$-combination of the polar syzygies ${\bf p}_{ijk}$ and of the polar syzygies that generate
$\widetilde{\fz(m-1)}$.
This shows the contention.
\qed

\medskip

We next treat a well-known class of scrollar  parametrizations.
The approach of \cite{Conca} is suited here, whereby any $d$-catalecticant $2\times m$ generic matrix
can be written as a scrollar matrix. As remarked in [loc. cit.] this gives a subclass of all
scrollar matrices for which the corresponding subalgebras are Koszul algebras, hence defined
by quadratic relations (that include the appropriate Grassmann--Pl\"ucker relations).

Explicitly, fixing an integer $1\leq d\leq m$, let now
$$\XX(m,d):=\left(
\begin{array}{ccccccc}
x_1 & x_2 &  \cdots & x_{1+d} & \cdots & x_{m-1} & x_m  \\
x_{1+d} & x_{2+d} & \cdots & x_{1+2d} & \cdots & x_{m-1+d} & x_{m+d}
\end{array}
\right)
$$
Note that the extreme values $d=1$ and $d=m$ correspond, respectively, to the
ordinary Hankel matrix and the generic matrix.
As in the previous example, we let $\ff(m,d)$ denote the $2\times 2$ minors of $\XX(m,d)$
and consider the $k$-subalgebra $A=k[\ff(m,d)]\subset R=k[x_1,\ldots,x_{m+d}]$.
It has also been shown that $\dim A=\min\{2m-3,m+d\}$ (see \cite[Theorem 5.2]{Conca}).
Therefore, $\dim A=\dim R$ if and only if $d\leq m-3$; for $d>m-3$, $A$ is isomorphic to the
homogeneous coordinate ring of the Grassmannian $\mathbb{G}(1,m+d-1)$ dicussed in the previous
example.

\begin{Proposition}\label{scrolls}
Let $\ff(m,d)$ denote the $2\times 2$ minors of the $d$-catalecticant $2\times m$ generic matrix $\XX(m,d)$.
Then $\ff(m,d)$ is polarizable.
\end{Proposition}
\demo
To see that $\ff(m,d)$ is polarizable, we would like to apply the same recipe as in the generic case.
This is fine insofar that removing the last column of the above matrix preserves its catalecticant shape with the
same leap $d$, thus allowing to induct
on the number of columns as in the completely generic case.

However, there are some numerical changes when $d\leq m-4$.
Firstly, the ideal $P(m,d)$ is now generated by the Grassmann--Pl\"ucker relations and additional
quadratic relations containing more than $3$ terms. According to \cite[Corollary 3.4]{Conca}
these latter relations -- called $N$-{\em relations} -- constitute a set whose cardinality is the number $\nu(m,d)$
of indices $i,j,k,l$ such that $1\leq i<j<k<l\leq m$ and $k-j>d$. Therefore, $\nu(m,d)={{m-d}\choose{4}}$,
and hence $\mu(P(m,d))={{m}\choose{4}} + {{m-d}\choose{4}}$.

A set of generators of $P(m,d)$ as given in \cite[Corollary 3.4]{Conca} can be further slightly modified so that
each $N$-relation is replaced by a quadratic relation with exactly $6$ terms, which for convenience we will call an
$M$-{\em relation}.

The advantage of switching to $M$-relations stems from the fact that, while an $N$-relation as
in \cite{Conca} comes naturally from a parallel relation at the level of the initial ideal, an $M$-relation is
seen to be structurally defined. Namely, for every choice of indices $i,j,k,l$ such that $d+1\leq i<j<k<l\leq m$,
pick the $2\times 4$ submatrix of $\XX(m,d)$ whose entries on the {\em first} row are $x_i,x_j,x_k,x_l$.
(Clearly, the total number of such different choices is $\nu(m,d)={{m-d}\choose{4}}$.)
Then pick the unique  $2\times 4$ submatrix of $\XX(m,d)$ whose entries on the {\em second} row are $x_i,x_j,x_k,x_l$.
Next stack these two matrices together to form a $4\times 4$ matrix $M(i,j,k,l)$ with one repeated row, so its determinant vanishes.
Expanding this determinant by the  Laplacian rule along the $2$-minors of the first two rows yields a quadratic relation between
a subset of the $2$-minors of $M(i,j,k,l)$ which, by construction, gives a relation between a subset of the $2$-minors
of the original matrix $\XX(m,d)$.
In this way, one gets exactly one $M$-relation for each choice of a $2\times 4$
submatrix of $\XX(m,d)$ involving four of the last $m-d$ columns of $\XX(m,d)$.

Now one argues in a similar way as in the generic case by a slight adaptation of the procedure.
\qed

\subsection{A class of polar maps}

This example is inspired by a renowned
 construction of Gordan and Noether
 in connection with the Hesse problem (see \cite{GN}).
 We will actually look only at the simplest situation of their
 construction as follows: let $f_1,\ldots,f_{n-r}\, (n\geq r+1)$ be forms of the same degree $\geq 2$ in
the polynomial ring $R:=k[x_1,\ldots,x_r]$ over a field $k$ -- in \cite{GN} the authors assume that the given forms are
algebraically dependent over $k$, e.g., when $n\geq 2r+1$, but will make no such restriction at the outset.

Consider the $k$-subalgebra $A:=k[f_1,\ldots,f_{n-r}]\subset R$.
Take  new variables $x_{r+1},\ldots,x_n$, and let $F$ be the
following form
$$
F:=x_{r+1}f_1+\ldots+x_nf_{n-r}\in S:=k[x_1,\ldots,x_r,x_{r+1},\ldots,x_n],
$$
which can suggestively  be thought of as the generic member of the linear system spanned by
$\ff:=\{f_1,\ldots,f_{n-r}\}$.

Let $B$ be the $k$-subalgebra of $S$ generated by the partial
derivatives of $F$.
Write
$$B:=k\left[\frac{\partial F}{\partial x_1},\ldots,\frac{\partial F}{\partial x_r},\frac{\partial F}{\partial x_{r+1}}
=f_1,\ldots,\frac{\partial F}{\partial x_{n}}=f_{n-r}\right]\subset S.$$
For simplicity, call these generators $\mathbf{g}=\{g_1,\ldots,g_r,g_{r+1},\ldots,g_n\}$.

\begin{Proposition}\label{gordan_noether}
With the above notation, if $\dim A=r$ then $B$ is polarizable if (and only if) $A$ is polarizable.
\end{Proposition}
\demo
Note that the hypothesis of the claim implies that $n\geq 2r$.

Actually, we prove something more precise.
Namely, look at the transposed Jacobian matrix of $\mathbf{g}$ which is nothing but the
Hessian matrix of $F$. A careful inspection of $\mathbf{g}$ leads to the following block shape

\arraycolsep=10pt

$$^t \Theta(\mathbf{g})=
\left(
\begin{array}{c@{\quad\vrule\quad}c}
\raise5pt\hbox{$\mathcal{H}$} & \raise5pt\hbox{$^t\Theta(\ff)$}\\[-6pt]
\multispan2\hrulefill\\
\Theta(\ff)&0
\end{array}
\right),
$$
where $\Theta(\ff)$ is the Jacobian matrix of the originally given forms $\ff$
and $\mathcal{H}$ is the $r\times r$ Hessian matrix of $F$ regarded as a polynomial in $x_1,\ldots,x_r$
with coefficients in $k[x_{r+1},\ldots, x_n]$.
Since $\dim A=r\leq n-r$ then the rank of $\Theta(\ff)$ is $r$, hence one readily obtains that
the first $r$ coordinates of any syzygy of the above matrix are all zero, while the last
$n-r$ coordinates constitute a syzygy of $^t\Theta(\ff)$.

To translate this outcome in a more formal fashion, let $\fz_{_T}(\mathbf{h})$ stand for the differential syzygy module
of a set of forms $\mathbf{h}$
in a polynomial ring $T$ over $k$.
By definition, we have the exact sequence of $S$-modules
\begin{equation}\label{transposed_jac}
0\rar \fz_{_R}(\ff)\lar R^{n-r}\stackrel{^t\Theta(\ff)}{\lar} \sum_{i=1}^r R\,dx_i.
\end{equation}
Since $R\subset S$ is a flat (free) extension, we get an exact sequence of $R$-modules
$$0\rar \fz_{_R}(\ff)\otimes _R S\lar S^{n-r}\stackrel{^t\Theta(\ff)\otimes_R S}{\lar} \sum_{i=1}^r S\,dx_i.$$

Taking $R$-duals in (\ref{transposed_jac}), observing that coker$(^t\Theta(\ff))$ is torsion because the rank
of $^t\Theta(\ff)$ is $r$, and tensoring with $S$ over $R$, yields a short exact complex
\begin{equation}\label{jac}
0\rar \sum_{i=1}^r S\,\frac{\partial}{\partial x_i} \stackrel{\Theta(\ff)\otimes_R S}{\lar} (S^{n-r})^*.
\end{equation}
Consider the map of complexes over $S$
{\small
$$\begin{array}{ccccccccc}
0 \kern-8pt & \kern-8pt \rar\kern-8pt & \kern-8pt\fz_{_R}(\ff)\otimes _R S\kern-8pt & \lar & S^{n-r} \kern-8pt
&\stackrel{^t\Theta(\ff)\otimes_R S}{\lar} \kern-8pt & \sum_{i=1}^r S\,dx_i\kern-8pt &&\\[3pt]
&&&& \kern-8pt\uparrow \kern-8pt && \uparrow\mathcal{H} \kern-8pt && \\[3pt]
&&&&    0 & \lar & \sum_{i=1}^r S\,\frac{\partial}{\partial x_i} & \stackrel{\Theta(\ff)\otimes_R S}{\lar}  &(S^{n-r})^*
\end{array}
$$
}
The mapping cone of this map is the exact complex
$$0\rar \fz_{_R}(\ff)\otimes _R S \lar S^{n-r}\oplus \sum_{i=1}^r S\,\frac{\partial}{\partial x_i}
 \stackrel{\Psi}{\lar} \sum_{i=1}^r S\,dx_i\oplus (S^{n-r})^*,$$
 where $\Psi$ is represented, in suitable bases, by the matrix

 \arraycolsep=10pt

$$
\left(
\begin{array}{c@{\quad\vrule\quad}c}
\raise5pt\hbox{$^t\Theta(\ff)$} & \raise5pt\hbox{$\mathcal{H}$}\\[-6pt]
\multispan2\hrulefill\\
0 & \Theta(\ff)
\end{array}
\right).
$$
Therefore, up to change of basis in the free modules, one gets an isomorphism
$\fz_{_R}(\ff)\otimes _R S\simeq \fz_{_S}(\mathbf{g})\subset S^n$.

In particular, these two $S$-modules have the same rank.
Now, the $S$-rank of $\fz_{_R}(\ff)\otimes _R S$ equals
the $R$-rank of $\fz_{_R}(\ff)$ and the latter is $n-r-\dim A$.
But the rank of the $S$-module $\fz_{_S}(\mathbf{g})$ is $n-\dim B$.
It follows that $\dim B=\dim A+r$. But since $B$ is generated over $A$
by $r$ elements, it follows that $B$ is a polynomial ring over $A$
(of dimension $2r$).

We now deal with the respective modules of polar syzygies.
Let $A\simeq k[T_1,\ldots,T_{n-r}]/P$ and $B\simeq k[T_1,\ldots,T_n]/\boldsymbol P$ denote
respective polynomial presentations of the two algebras.
One has an inclusion $Pk[T_1,\ldots,T_n]\subset\boldsymbol P$.
Since $\hht (P)=n-r-\dim A=n-2r=n-\dim B=\hht(\boldsymbol P)$, we have an equality.
This implies a similar relation as above between the respective associated modules of polar syzygies.
Namely,  $B$ is a polynomial ring over $A$, hence $A\subset B$ is a free extension.
Therefore, by the same token, the equality $Pk[T_1,\ldots,T_n]=\boldsymbol P$ translates into
a natural isomorphism
$${\rm Im}(P/P^2)\otimes_A B\simeq {\rm Im}(\boldsymbol P/{\boldsymbol P}^2)\subset \sum_{\ell=1}B\, dT_{\ell}.$$
This shows the contention.
\qed

\bigskip

\noindent {\sc Illustration.}
Take $A$ to be a $d$-catalecticant parameterization as in (\ref{scrolls}), with $d\leq m-3$.
Then, as remarked there, $A$ has maximal dimension, i.e., $\dim A=m+d$.
The polar map of the corresponding hypersurface $F$ has an image of (geometric) dimension
$2m+2d-1$ in $\pp^{{{m}\choose {2}} +m+d-1}$.
As a slight check on the numbers, since we are assuming $d\leq m-3$ then $m+d\leq 2m-3$. Clearly
then ${{m}\choose {2}}< m+d\leq 2m-3$ would entail $(m-2)(m-3)=m^2-5m+6<0$, which is impossible
as $m\geq 3$.

In the simplest case, with $d=1$ and $m=4$,we get the polarizable forms
{\small
$$
\begin{array}{l}
x_3x_6+x_4x_7+x_5x_8\boldsymbol, \; -2x_2x_6-x_3x_7-x_4x_8+x_4x_9+x_5x_{10}\boldsymbol, \;\; x_1x_6-x_2x_7-2x_3x_9-\\
-x_4x_{10}+x_5x_{11}\boldsymbol,\;\; x_1x_7-x_2x_8+x_2x_9-x_3x_{10}-2x_4x_{11}\boldsymbol, \;\;  x_1x_8+x_2x_{10}+x_3x_{11}\boldsymbol, \\
-x_2^2+x_1x_3\boldsymbol,\; -x_2x_3+x_1x_4\boldsymbol,\; -x_2x_4+x_1x_5\boldsymbol, \; -x_3^2+x_2x_4\boldsymbol,\;
 -x_3x_4+x_2x_5\boldsymbol,\; -x_4^2+x_3x_5
\end{array}
$$
}
parameterizing a Pl\"ucker quadric hypersurface in $\pp^{10}$.

\begin{Remark}\rm We observe that this procedure is in principle iterative, but the numbers and
the result proper will be different.
\end{Remark}

We close with some general questions.

\begin{Question}
Let $\ff\subset R$ be a polarizable set of forms of degree $2$.
Is $k[\ff]\subset R$ normal?
\end{Question}

\begin{Question}
Let $\ff\subset R$ be a polarizable set of forms of degree $2$.
Is $k[\ff]\subset R$ Cohen--Macaulay?
\end{Question}

In \cite{BGS} the answer to the first of these questions is affirmative when
$\ff$ are monomials, hence so is
the answer to the second question under the same assumption.

In general, if the answer to the second of these questions is affirmative -- or, if at least
$k[\ff]$ satisfies Serre's property $(S_2)$ -- a loose strategy for the first question would
be to show that $\hht \mathfrak{J}R\leq \hht \mathfrak{J}$. This would entail $\hht \mathfrak{J}\geq 2$
since with $\fp=\fz$ the ideal $\mathfrak{J}R$ would be the Fitting of a second syzygy.
Thus, a preliminary question along this line is to understand when the Jacobian ideal
is height decreasing through $k[\ff]\subset R$.
A sufficient condition for this to happen is that the contraction $\mathfrak{J}R\cap k[\ff]$
be contained in some minimal prime of $\mathfrak{J}$ of minimal height.

Affirmative answers to these questions would give another proof of
the result of \cite{Conca} to the effect that the algebra parameterized by
the $2\times 2$ minors of a $d$-catalecticant matrix is normal and Cohen--Macaulay.
Of course, this is not a totally impressive outcome since we have drawn upon the ideas of [loc.cit.] for the structure
of such objects in order to prove their polarizability.

\begin{Question}\label{koszul}\rm
Suppose  that $P$ is generated by forms of degree $2$
constituting a Gr\"obner basis for some monomial order. Is $\ff$ polarizable?
\end{Question}
The assumption forces the $k$-algebra $A=k[\ff]\subset R$ to be a Koszul algebra
according to \cite[Theorem 2.2]{BHV}.
Note that most geometric examples in this section are Koszul algebras.


\begin{thebibliography}{99}

\bibitem{BGS}{I. Bermejo, P. Gimenez and A. Simis, Polar syzygies
 in characteristic zero: the monomial case, J. Pure Appl. Algebra
{\bf 213} (2009) 1--21.}

\bibitem{gauss}{P. Brumatti, P. Gimenez and A. Simis, On the Gauss algebra associated to a rational
map $\pp^d \dasharrow \pp^n$, J. Algebra {\bf 207} (1998), 557--571.}

\bibitem{BHV}{W. Bruns, J. Herzog and U. Vetter, Syzygies and walks, {\em in} ICTP Proceedings ``Commutative
Algebra'', Eds. A. Simis, N. V. Trung and G. T. Valla, World Scientific (1994) 36--57.}

\bibitem{Conca}{A. Conca, J. Herzog and G. Valla, Sagbi bases with applications to blow-up algebras,
J. reine angew. Math. {\bf 474} (1996), 113--138.}

\bibitem{GN}{P. Gordan und M. Noether, Ueber die algebraischen Formen, deren
Hesse'sche Determinante identisch verschwindet, Math. Ann., {\bf 10}
(1876), 547--568.}


\bibitem{He78}{J. Herzog, Ein Cohen-Macaulay Kriterium mit
Anwendungen auf den Konormalenmodul und den Differentialmodul, Math.
Z. {\bf 163} (1978), 149--162. }


\bibitem{jac}{A. Simis, On the jacobian module associated to a graph,
Proc. Amer. Math. Soc., {\bf 126} (1998), 989--997.}

\bibitem{elam99}{A. Simis, {\it Remarkable Graded Algebras in Algebraic Geometry},
{\sc XII ELAM}, IMCA, Lima, Peru, 1999.}

\bibitem{Mich}{A.  Simis, Two differential themes in characteristic zero, {\it in\/} {\sc Topics
in Algebraic and Noncommutative Geometry}, Proceedings in Memory of Ruth Michler
(Eds. C. Melles, J.-P. Brasselet, G. Kennedy, K. Lauter and L. McEwan),
Contemporary Mathematics {\bf 324}, Amer. Math. Soc., Providence, RI, 2003, 195--204.}

\bibitem{ram3}{A. Simis, B. Ulrich and W. Vasconcelos, Tangent algebras,
Trans. Amer. Math. Soc., {\bf 364} (2012), 571--594.}



\end{thebibliography}
\end{document}